\documentclass[11pt, reqno]{amsart}

\usepackage{amsthm, amsmath, amssymb,xcolor,graphicx}

\usepackage[hidelinks]{hyperref}
\usepackage{url}

\newtheorem{theorem}{Theorem}

\newtheorem{conjecture}[theorem]{Conjecture}
\newtheorem*{question*}{Question}

\theoremstyle{definition}

\newtheorem*{definition*}{Definition}

\theoremstyle{remark}

\DeclareMathOperator{\Aff}{Aff}

\DeclareMathOperator{\SL}{SL}

\newcommand{\FF}{\mathbb{F}}

\newcommand{\PP}{\mathbb{P}}

\newcommand{\ZZ}{\mathbb{Z}}

\title{Group representations that resist worst-case sampling}
\author{Yufei Zhao}
\address{Mathematical Institute, Oxford OX2 6GG, United Kingdom}
\email{yufei.zhao@maths.ox.ac.uk}
\thanks{The author was supported by an Esm\'ee Fairbairn 
Junior Research Fellowship at New College, Oxford and a Research Fellowship at the Simons Institute, Berkeley.}

\begin{document}

\begin{abstract}
Motivated by expansion in Cayley graphs, we show that there exist infinitely many groups $G$ with a nontrivial irreducible unitary representation whose average over every set of $o(\log\log|G|)$ elements of $G$ has operator norm $1 - o(1)$. This answers a question of Lovett, Moore, and Russell, and strengthens their negative answer to a question of Wigderson. 

The construction is the affine group of $\mathbb{F}_p$ and uses the fact that for every $A \subset \mathbb{F}_p\setminus\{0\}$, there is a set of size $\exp(\exp(O(|A|)))$ that is almost invariant under both additive and multiplicatpive translations by elements of $A$.
\end{abstract}

\maketitle

Let $G$ be a finite group and $\rho$ a unitary representation of $G$. For a subset $S \subset G$, we say that $S$ is \emph{$\epsilon$-expanding} with respect to $\rho$ if
\[
\biggl\| \frac{1}{|S|}\sum_{g \in S} \rho(g) \biggr\|_{\mathrm{op}} \le 1 - \epsilon.
\]
Otherwise, we say that $\rho$ \emph{$\epsilon$-resists} $S$. We say that $S \subset G$ is \emph{$\epsilon$-expanding} if it is $\epsilon$-expanding with respect to every non-trivial irreducible unitary representation of $G$, which is essentially the same as saying that the adjacency matrix of the Cayley graph on $G$ generated by $S$ has all eigenvalues, except the top one, bounded by $1-\epsilon$ in absolute value. It is closely related to a more combinatorial notion of expansion in graphs via Cheeger's inequality.

By a theorem of Alon and Roichman~\cite{AR94} on eigenvalues of random Cayley graphs, for any group $G$, a random set of $C_\epsilon \log |G|$ group elements is $\epsilon$-expanding with high probability. This bound is tight for abelian groups, up to constant factors. For example, when $G = (\ZZ/2\ZZ)^n$, it takes $n = \log_2 |G|$ elements simply to generate the group.

On the other hand, for certain families of ``highly non-abelian'' groups, including all non-abelian simple groups, a bounded number of generators suffices to obtain $\epsilon$-expansion.  In certain cases, such as $\SL_2(\FF_p)$~\cite{BG08}, and more generally, any finite simple groups of Lie type of bounded rank~\cite{BGGT15}, we know that $\{g,h,g^{-1},h^{-1}\}$ is $\epsilon$-expanding with high probability for uniformly random group elements $g$ and $h$. See surveys~\cite{HLW06,Lub12} for more on expansion.

Wigderson conjectured~\cite{Wig} in his 2010 Barbados lectures that there is some constant $k$ so that for any finite group $G$ and a nontrivial irreducible unitary representation $\rho$, a list of $k$ random elements of $G$ is $(1/2)$-expanding with respect to $\rho$ with probability at least $1/2$. Note that this is true for abelian groups, where every irreducible representation is one-dimensional, even though it takes $C \log |G|$ elements to expand with respect to every non-trivial irreducible representation simultaneously.

Wigderson's conjecture was disproved by Lovett, Moore, and Russell~\cite{LMR15}, who found an infinite family of groups such that, with high probability, a random subset of $k$ elements does not expand at all with respect to a specific nontrivial irreducible representation. More specifically, they showed that if $K$ is a fixed non-abelian group with trivial center (e.g., $K = S_3$), and $\rho$ is a faithful irreducible unitary representation of $\rho$, then, $G = K^n$ with the irreducible representation $\rho^n = \rho \otimes \cdots \otimes \rho$ has the property that, as $n \to \infty$, provided that $k = o(\log n)$, one has
\[
\PP_{g_1, \dots, g_k \in G} \left[ \left\| \frac{\rho^n(g_1) + \cdots + \rho^n(g_k)}{k} \right\|_{\mathrm{op}} = 1\right] = 1 - o(1).
\]
Therefore, there are infinitely many groups $G$ with a non-trivial irreducible unitary representation that resist any set of size $o(\log\log|G|)$. Despite these negative results for random group elements, they asked whether there are constants $k$ and $\epsilon$ such that for any group $G$  and any nontrivial irreducible representation $\rho$, there exist some $k$ elements of $G$ that $\epsilon$-expand respect to $\rho$. We answer this question in the negative. 

\begin{theorem} \label{thm:resistant}
For every $\epsilon >0$, there is some $c_\epsilon > 0$ so that there exist infinitely many groups $G$ with a nontrivial irreducible unitary representation $\rho$ that $\epsilon$-resists every $S\subset G$ with $|S| \le c_\epsilon \log\log|G|$, i.e., 
\begin{equation} \label{eq:norm}
\left\| \frac{\rho(g_1) + \cdots + \rho(g_k)}{k} \right\|_{\mathrm{op}} \ge 1 - \epsilon
\end{equation}
for any $g_1, \dots, g_k \in G$ with $k \le c_\epsilon \log \log |G|$.
\end{theorem}

More succinctly, there exist groups $G$ with a representation that $o(1)$-resists any set of $o(\log\log|G|)$ elements (the construction in \cite{LMR15} works for a random set, whereas ours works for all sets).

This gives a strong negative answer to Wigderson's question, as it shows that there no choice of a constant number of elements of $G$ can $\epsilon$-expand with respect to $\rho$, let alone a random choice.

\medskip

We prove Theorem~\ref{thm:resistant} by taking $G = \Aff(\FF_p)$, the affine group of $\FF_p$. Its elements are affine transformations $x \mapsto ax + b$, where $a \in \FF_p^\times$ and $b \in \FF_p$. Let $\rho$ denote its standard representation with the trivial component removed. Theorem~\ref{thm:resistant} for $\Aff(\FF_p)$ is an immediate consequence of the following result.

\begin{theorem}
	\label{thm:affine}
	For every $\epsilon > 0$, there is some $C_\epsilon > 0$ so that for every prime $p$ and every $A \subset \FF_p\setminus\{0\}$, there exists some $X \subset \FF_p$ with $|X| \le \exp(\exp(C_\epsilon |A|))$ such that
	\begin{equation}
		 \label{eq:affine-invar}
		|a \cdot X \setminus X| \le \epsilon |X|
		\quad \text{and} \quad
		|(a + X) \setminus X| \le \epsilon |X|
		\quad \text{for all } a \in A.
	\end{equation}
\end{theorem}
Here we adopt the standard notation from additive combinatorics: $a \cdot X := \{ax : x \in X\}$ and $a + X := \{a +x : x \in X\}$. Theorem~\ref{thm:resistant} for $G = \Aff(\FF_p)$ follows as a corollary of Theorem~\ref{thm:affine}. Indeed, if $S$ consists of affine maps $x \mapsto a_i x + b_i$, then take $A$ to be the set of all nonzero elements that appears as $a_i$ or $b_i$ for some $i$. The claim~\eqref{eq:norm} follows by considering the characteristic vector of $X$, appropriately normalized (noting $|X| \le p/2$ if $|S| \le c_\epsilon \log\log p$; we may need to rescale $\epsilon$ by a constant factor).

A proof of Theorem~\ref{thm:affine} was given by Terry Tao in a MathOverflow post~\cite{TaoMO}.\footnote{The author thanks Ben Green for pointing out \cite{TaoMO} to him.} We include the proof here for completeness.

\begin{proof}
	Let $A = \{a_1, \dots, a_k\}$.
	Let $L = \lceil 1/\epsilon \rceil$. Consider the generalized arithmetic and geometric progressions
	\begin{align*}
		P &= \{ n_1 a_1 + \cdots + n_k a_k : 0 \le n_1, \dots, n_k < L\}, \text{ and}\\ 
		Q &= \{ a_1^{n_1} \cdots a_k^{n_k} :  0 \le n_1, \dots, n_k < L\}.
	\end{align*}
	Let
	\[
	X = Q^{-1} \biggl( \sum_{y \in Q} y \cdot P \biggr),
	\]
	i.e., the set of all elements that can be written as
	\[
	y_0^{-1} \biggl( \sum_{y \in Q} y x_y \biggr)
	\]
	for some choices of $y_0 \in Q$ and $x_y \in P$ for each $y \in Q$. It is easy to check \eqref{eq:affine-invar}, as $|(a + P)\setminus P| \le \epsilon |P|$ and $|(a \cdot Q)\setminus Q| \le \epsilon |Q|$ for any $a \in A$. We have $|X| \le |Q| |P|^{|Q|} \le L^k L^{kL^k} \le e^{(1/\epsilon)^{O(k)}}$.
\end{proof}

It remains an open question whether the bounds in Theorems~\ref{thm:resistant} and \ref{thm:affine} can be improved. We conjecture that they cannot.

\begin{conjecture}  \label{conj:rep}
	For every $\epsilon > 0$, there is some $C_\epsilon$ such that for any group $G$ and a nontrivial irreducible unitary representation $\rho$, there is some $S \subset G$ with $|S| \le C_\epsilon \log\log|G|$ that is $\epsilon$-expanding with respect to $\rho$.
\end{conjecture}

\begin{conjecture}
	For every $\epsilon > 0$, there is some $c_\epsilon > 0$ so that for every positive integer $k \le c_\epsilon \log \log p$ and prime $p$, there is some $A \subset \FF_p\setminus \{0\}$ with $|A| = k$ such that every nonempty $X \subset \FF_p$ satisfying \eqref{eq:affine-invar} has $|X| \ge \exp(\exp(c_\epsilon k))$.
\end{conjecture}

\begin{conjecture}
	In the above conjectures, choosing $S$ and $A$ uniformly at random works with high probability.
\end{conjecture}

Note that Alon--Roichman theorem implies that Conjecture~\ref{conj:rep} is true if we replace $\log\log|G|$ by $\log |G|$ (by taking a random $S$).
 We do not know any further improvements.

\subsection*{Acknowledgment} The author thanks Ben Green for discussion and Shachar Lovett for encouraging him to write up this result.

%\bibliographystyle{amsplain_mod2}
%\bibliography{ref_resist}

\end{document}